\documentclass[10pt]{article}
\usepackage{amsmath,amsthm,amsfonts,amssymb}
\usepackage{epsfig}
\usepackage{color}
\usepackage{latexsym}
\usepackage{supertabular}
\usepackage{graphicx}
\usepackage{a4wide}
\numberwithin{equation}{section}

\def\whitebox{{\hbox{\hskip 1pt
 \vrule height 6pt depth 1.5pt
 \lower 1.5pt\vbox to 7.5pt{\hrule width
    3.2pt\vfill\hrule width 3.2pt}%
 \vrule height 6pt depth 1.5pt
 \hskip 1pt } }}
\def\qed{\ifhmode\allowbreak\else\nobreak\fi\hfill\quad\nobreak
     \whitebox\medbreak}

\newcommand{\ignore}[1]{}

\theoremstyle{plain}
\newtheorem{theorem}{Theorem}[section]
\newtheorem{corollary}[theorem]{Corollary}
\newtheorem{lemma}[theorem]{Lemma}

\newtheorem{problem}[theorem]{Problem}
\newtheorem{proposition}[theorem]{Proposition}
\newtheorem{example}[theorem]{Example}

\newtheorem{remark}[theorem]{Remark}

\def\qed{{\hfill$\square$}}
\def\proof{{\vspace{-0.3cm}\bf Proof: \,}}

\def\Z{{\mathbb Z}}
\def\Q{{\mathbb Q}}

\def\C{{\mathbb C}}
\def\F{{\mathbb F}}
\def\mod{{\mathrm{mod\,\,}}}

\def\Tr{{\mathrm{Tr}}}

\def\Cay{{\mathrm{Cay}}}

\title{Skew Hadamard difference sets from cyclotomic strongly regular graphs}
\author{Koji Momihara\footnotemark[1]}
\date{} 
\begin{document}
\maketitle
\footnotetext[1]{
Department of Mathematics, Faculty of Education, Kumamoto University,  
2-40-1 Kurokami, Kumamoto 860-8555, Japan; Email address: 
momihara@educ.kumamoto-u.ac.jp}
\renewcommand{\thefootnote}{\arabic{footnote}}
\begin{abstract}  
We find new constructions of infinite families of skew Hadamard difference sets in elementary abelian groups under the assumption of  the existence of cyclotomic strongly regular graphs. 
Our construction is based on choosing
cyclotomic classes in finite fields. 
\end{abstract}
\begin{center} 
\end{center}
\section{Introduction}
We assume that the reader is familiar with the basic theories of  difference sets and strongly regular graphs (srg) as can
found in \cite{bjl,bh}. 

A difference set $D$ in an (additively written) finite group $G$ is 
called {\it skew Hadamard} if $G$ is the disjoint union of $D$, $-D$, and $\{0\}$. 
The primary example (and for many years, the only known example in abelian groups) of 
skew Hadamard difference sets is the classical Paley (quadratic residue) 
difference set in $(\F_q,+)$ consisting of the nonzero squares of $\F_q$, 
where $\F_q$ is the finite field of order $q$, a prime power congruent to 3 modulo 4.  
Skew Hadamard difference sets are currently under intensive study, see  
\cite{CP,DY06,DWX07,F11,FX113,FMX11,M,WH09,WQWX}. There were two major 
conjectures in this area: (i) If an abelian group $G$
contains a skew Hadamard difference set, then $G$ is necessarily elementary
abelian. (ii) Up to equivalence the Paley difference sets mentioned above
are the only skew Hadamard difference sets in abelian groups. The former
conjecture is still open in general. 
The latter conjecture turned out to be false: Ding
and Yuan \cite{DY06} constructed a family of skew Hadamard difference sets in
$(\F_{3^m},+)$, where $m\ge 3$ is odd, and showed that two examples
in the family are inequivalent to the Paley difference sets. Very recently,
Muzychuk \cite{M} constructed infinitely many inequivalent skew Hadamard
difference sets in an elementary abelian group of order $q^3$. The reader may check the introduction of \cite{FX113} for a good short survey of known constructions of skew Hadamard difference sets and related problems.

A classical method for constructing both connection sets of strongly regular Cayley graphs (called {\it partial difference sets}) and ordinary difference sets in the additive groups of finite fields is to use cyclotomic classes of finite fields. 
Let $p$ be a prime,  $f$ a positive integer, and let $q=p^f$. Let $k>1$ be an integer such that $k|(q-1)$, and $\gamma$ be a primitive root of $\F_q$. 
Then the cosets $C_i^{(k,q)}=\gamma^i \langle \gamma^k\rangle$, $0\leq i\leq k-1$, are called the {\it cyclotomic classes of order $k$} of $\F_q$. 
Many authors have studied the problem of determining when a union $D$ of some cyclotomic classes forms a (partial) difference set.
Especially, when $D$ consists of only a subgroup of $\F_q$, 
many authors have studied extensively \cite{BM73,BMW82,FX111,FX113,FMX11,GXY11,M75,SW02,S67,VLSch}. 
We call such a strongly regular Cayley graph $\Cay(\F_q,D)$ {\it cyclotomic}. 
The well known Paley graphs are primary examples of  cyclotomic srgs.  
Also, if $D$ is the multiplicative group of a
subfield of $\F_q$, then it is clear that $\Cay(\F_q , D)$ is strongly regular. 
These cyclotomic srgs are usually called {\it subfield examples}. Next, if there exists a
positive integer $t$ such that $p^t\equiv -1\,(\mod{k})$, then $\Cay(\F_q , D)$ is strongly regular. 
This case is usually called {\it semi-primitive}. In \cite{SW02}, 
Schmidt and White conjectured that if  $k\,|\,\frac{q-1}{p-1}$ and $\Cay(\F_{p^f},C_0^{(k,q)})$ is strongly regular, then one of the following holds: 
\begin{enumerate}
\item[(1)] (subfield case) $C_0=\F_{p^d}^\ast$ where $d\,|\,f$,  
\item[(2)] (semi-primitive case) $-1\in \langle p\rangle\le (\Z/k\Z)^\ast$,
\item[(3)] (exceptional case) $\Cay(\F_{p^f},C_0)$ has one of the parameters 
given in Table~\ref{Tab1}. 
\begin{table}[h]
\caption{Eleven sporadic examples}
\label{Tab1}
$$
\begin{array}{|c||c|c|c|c|}
\hline
\mbox{No.}&k&p&f&e:=[(\Z/k\Z)^\ast:\langle p\rangle]\\
\hline
1&11&3&5&2\\
2&19&5&9&2\\
3&35&3&12&2\\
4&37&7&9&4\\
5&43&11&7&6\\
6&67&17&33&2\\
7&107&3&53&2\\
8&133&5&18&6\\
9&163&41&81&2\\
10&323&3&144&2\\
11&499&5&249&2\\
\hline
\end{array}
$$
\end{table}
\end{enumerate}
Recently, in \cite{FX111,FMX11,GXY11,M12}, it was succeeded to generalize the sporadic examples of Table~\ref{Tab1} except for the srg of No.~1 and several subfield examples into 
infinite families using ``index $2$ or $4$ Gauss sums'' and ``relative Gauss sums.'' Also, Wu \cite{Wu12} gave a necessary and sufficient condition for 
$\Cay(\F_{p^{(k-1)/e}},C_0^{(p_1,p^{(k-1)/e})})$ to be strongly regular by generalizing the method of \cite{GXY11} when $k$ is a prime. 
On the other hand, in \cite{FX113,FMX11}, Feng, Xiang, and this author found new constructions of 
skew Hadamard difference sets via a computation of a character sum involving index $2$ Gauss sums. In particular, in \cite{FX113,FMX11}, it was shown  that $D=\bigcup_{i\in \{0\}\cup \langle p\rangle\cup 2\langle p\rangle}C_i^{(k,p^f)} $ is a skew Hadamard difference sets or a Paley type partial difference sets for the triples $(k,p,f)$ of Table~\ref{Subsec2Tab1} and these examples can be generalized into infinite families.  (A partial difference set $D$ in a group $G$ is said
to be of {\it Paley type} if the parameters of the corresponding strongly regular Cayley graph are
$(v, (v-1)/2,(v-5)/4,(v-1)/4)$.)  
\begin{table}[h]
\begin{center}
\caption{\label{Subsec2Tab1}
Skew Hadamard difference sets and Paley type partial difference sets  from index $2$ case}
$$
\begin{array}{|c||c|c|c|}
\hline
\mbox{No.} &\mbox{$k$}&\mbox{$p$} & \mbox{$f$}\\
\hline \hline 
\mbox{1} & \mbox{$2\cdot 11$}&\mbox{$3$} & \mbox{$5$}\\
\mbox{2} & \mbox{$2\cdot 19$}&\mbox{$5$} &\mbox{$9$}\\
\mbox{3} & \mbox{$2\cdot 67$}&\mbox{$17$} &\mbox{$33$}\\
\mbox{4} &\mbox{$2\cdot 107$}&\mbox{$3$}& \mbox{$53$}\\
\mbox{5} &\mbox{$2\cdot 163$}&\mbox{$41$}& \mbox{$81$}\\
\mbox{6} &\mbox{$2\cdot 499$}&\mbox{$5$}& \mbox{$249$} \\
\hline
\end{array}
$$
\end{center}
\end{table}
Now, one may recognize an interesting interaction between cyclotomic srgs and skew Hadamard difference sets of Tables~\ref{Tab1} and \ref{Subsec2Tab1}: for odd primes $p$ and $p_1$  such that $p$ is of index $2$ modulo $p_1$, the graph 
$\Cay(\F_{q},C_0^{(p_1,p^{(p_1-1)/2})})$ is 
strongly regular if and only if 
$D=\bigcup_{i\in \{0\}\cup \langle p\rangle\cup 2\langle p\rangle}C_i^{(2p_1,p^{(p_1-1)/2})} $ 
is  a skew Hadamard difference set or a Paley type partial difference set in 
$\F_q$. 

In this paper, we investigate such a relation between 
cyclotomic srgs and skew Hadamard difference sets, and 
find new constructions of infinite families of skew Hadamard difference sets  
from known cyclotomic srgs. 
\section{Background}
Let $p$ be a prime, $f$ a positive integer, and $q=p^f$. The canonical additive character $\psi$ of $\F_q$ is defined by 
$$\psi\colon\F_q\to \C^{\ast},\qquad\psi(x)=\zeta_p^{\Tr _{q/p}(x)},$$
where $\zeta_p={\rm exp}(\frac {2\pi i}{p})$ and $\Tr _{q/p}$ is the trace from $\F_q$ to $\F_p$. For a multiplicative character 
$\chi_k$ of order $k$ of $\F_q$, we define the {\it Gauss sum} 
\[
G_f(\chi_k)=\sum_{x\in \F_q^\ast}\chi_k(x)\psi(x), 
\] 
which belongs to the ring $\Z[\zeta_{kp}]$ of integers in the cyclotomic field $\Q(\zeta_{kp})$.
Let $\sigma_{a,b}$ be the automorphism of $\Q(\zeta_{kp})$ determined 
by 
\[
\sigma_{a,b}(\zeta_k)=\zeta_{k}^a, \qquad
\sigma_{a,b}(\zeta_p)=\zeta_{p}^b 
\]
for $\gcd{(a,k)}=\gcd{(b,p)}=1$. 
Below are several basic properties of Gauss sums \cite{BEW97}: 
\begin{itemize}
\item[(i)] $G_f(\chi_k)\overline{G_f(\chi_k)}=q$ if $\chi$ is nontrivial;
\item[(ii)] $G_f(\chi_k^p)=G_f(\chi_k)$, where $p$ is the characteristic of $\F_q$; 
\item[(iii)] $G_f(\chi_k^{-1})=\chi_k(-1)\overline{G_f(\chi_k)}$;
\item[(iv)] $G_f(\chi_k)=-1$ if $\chi_k$ is trivial;
\item[(v)] $\sigma_{a,b}(G_f(\chi_k))=\chi_k^{-a}(b)G_f(\chi_k^a)$.
\end{itemize}
In general, to  explicitly evaluate Gauss sums is very difficult. There are only a few cases where the Gauss sums have been evaluated. 
The most well known case is {\it quadratic} case, in other words, the order of $\chi$ is two. In this case, as can found in 
 \cite[Theorem 11.5.4]{BEW97}, it holds that 
\begin{equation}\label{eq:quad}
G_f(\chi_k)=(-1)^{f-1}\left(\sqrt{(-1)^{\frac{p-1}{2}}p}\right)^f. 
\end{equation}
The next simple case is the so-called {\it semi-primitive case} (also 
referred to as {\it uniform cyclotomy} or {\it pure Gauss sum}), where there 
exists an integer $j$ such that $p^j\equiv -1\,(\mod{k})$, where $k$ is the order of 
the multiplicative character $\chi$ involved. The explicit evaluation of Gauss sums in this case is given in \cite{BEW97}. 
The next interesting case is the index $2$ case where the subgroup $\langle p\rangle$ generated by $p\in (\Z/{k}\Z)^\ast$ is of index $2$ in $(\Z/{k}\Z)^\ast$ and $-1\not\in \langle p\rangle $. In this case, 
it is known that $k$ can have at most two odd prime divisors. 
Many authors have investigated this case, see \cite{YX10} for the complete solution  to the problem of evaluating index $2$ Gauss sums. 
Recently, these index $2$ Gauss sums were applied to 
show the existence of infinite families of 
new strongly regular graphs and skew Hadamard difference sets 
in $\F_q$ \cite{FX111,FX113,FMX11}. 

Now we recall the following well-known lemmas in the theories  of difference sets and strongly regular graphs
(see e.g., \cite{bh,M94}). 
\begin{lemma}\label{Sec3Le1}
Let $(G, +)$ be an abelian group of odd order $v$, $D$ a subset of $G$ 
of size $\frac{v-1}{2}$.  Assume that $D\cap -D=\emptyset$. Then, $D$ is a skew Hadamard difference set 
in $G$ if and only if 
\[
\psi(D)=\frac{-1\pm \sqrt{-v}}{2}
\] 
for all nontrivial characters $\psi$ of $G$.  On the other hand, assume that $0\not\in D$ and $-D=D$. Then $D$ is a Paley type partial difference set in $G$ if and only if 
\[
\psi(D)=\frac{-1\pm \sqrt{v}}{2}
\] 
for all nontrivial characters $\psi$ of $G$. 
\end{lemma}
\begin{lemma}\label{Sec3Le2}
Let $(G, +)$ be an abelian group and $D$ a subset of $G$.  Then, $\Cay(G,D)$ is a strongly regular graph  if and only if the size of the set 
\[
\{\psi(D)\,|\,\psi\in \widehat{G}\setminus \{\psi_0\}\}
\] 
is exactly two, where $\widehat{G}$ is the character group of $G$ and $\psi_0$ is the trivial character. 
\end{lemma}
Let $q$ be a prime power and let 
$C_i^{(k,q)}=\gamma^i \langle \gamma^k\rangle$, $0\le i\le k-1$, be 
the cyclotomic classes of order $k$ of $\F_q$, where $\gamma$ is a fixed primitive root of $\F_q$. In this paper, we assume that 
$D$ is a union of cyclotomic classes of order $k$ of $\F_q$. In order to check whether a candidate subset $D=\bigcup_{i\in I}C_i^{(k,q)}$ is a skew Hadamard difference set or a Paley type partial difference set, we will compute the sums $\psi(aD)=\sum_{x\in D}\psi(ax)$ for all $a\in \F_q^\ast$,  where $\psi$ is the canonical additive character of $\F_q$, because of Lemma~\ref{Sec3Le1}. Similarly, to check whether $D$ is a connection set of a strongly regular Cayley graph, we should compute the sums $\psi(aD)$ for all $a\in \F_q^\ast$ by Lemma~\ref{Sec3Le2}. 
Note that the sum $\psi(aD)$ can be expressed as a linear combination of 
Gauss sums  (cf. \cite[Lemma 3.1]{Wu12}) by using the orthogonality of characters: 
\begin{eqnarray*}
\psi(aD)=\frac{1}{k}
\sum_{\chi\in C_0^{\perp}}G_f(\chi^{-1})
\sum_{i\in I}\chi(a\gamma^i ), 
\end{eqnarray*}
where 
$C_0^{\perp}$ is the subgroup of $\widehat{\F_q^\ast}$
consisting of all $\chi$ which are trivial on $C_0^{(k,q)}$. 
Thus, the computation to know whether a candidate subset $D=\bigcup_{i\in I}C_i^{(k,q)}$ is a skew Hadamard difference set or a Paley type partial difference set is essentially  reduced to evaluating Gauss sums. 
In fact, in \cite{FX113,FMX11}, known evaluation of index $2$ Gauss sums are used. 
However, as previously said,  to  explicitly evaluate Gauss sums is very difficult. In this paper, we will show the existence of skew Hadamard difference sets and Paley type partial difference sets without computing explicit values of 
Gauss sums. 
Instead, we use the following theorem, called the 
{\it Davenport-Hasse product formula}   
\begin{theorem}
\label{thm:Stickel2}(\cite{BEW97})
Let $\eta$ be a multiplicative character of order $\ell>1$ of  $\F_q=\F_{p^f}$. For  every nontrivial character $\chi$ on $\F_q$, 
\[
G_f(\chi)=\frac{G_f(\chi^\ell)}{\chi^\ell(\ell)}
\prod_{i=1}^{\ell-1}
\frac{G_f(\eta^i)}{G_f(\chi\eta^i)}. 
\]
\end{theorem}
\section{Construction of skew Hadamard difference sets}
To show our main theorem, we use the following result of \cite{SW02}.  (They gave this result in terms of irreducible cyclic codes.)  
\begin{proposition}(\cite[Lemma~2.8, Corollary~3.2]{SW02})\label{prop}
Let $m$ be the order of $p$ modulo $k$ and set $q=p^f=p^{sm}$.  
Assume that $k\,|\,\frac{p^f-1}{p-1}$ is odd and  $\Cay(\F_{p^f},C_0^{(k,q)})$ is strongly regular. Then, for  a system $L$ of coset representatives of $\F_{p^f}^\ast/C_0^{(k,p^f)}$, there exists a partition 
$L_1\cup L_2=L$ such that  
\[
G_f(\chi_k)=\epsilon p^{s\theta}\sum_{x\in L_1}\chi_{k}(x)=
-\epsilon p^{s\theta}\sum_{x\in L_2}\chi_{k}(x), 
\]
where $\epsilon=\pm 1$ and $\theta$ is the integer such that 
$p^\theta|| G_{m}(\chi_k)$. 
(In this case, $p^{s\theta}||G_f(\chi_k)$ also holds.) Furthermore, if $|L_1|=k-d$ and 
$|L_2|=d$,  then it holds that 
\begin{eqnarray}\label{eq:Wu1}
k\cdot \psi(\gamma^a C_0^{(k,q)})+1&=&
\sum_{\chi\in C_0^{\perp\ast}}\chi(\gamma^a)G_f(\chi^{-1})\nonumber\\
&=&p^{s\theta}\epsilon d \mbox{\, \,  or\, \,  }
p^{s\theta}\epsilon (d-k). 
\end{eqnarray}
\end{proposition}
\begin{remark}{\em 
Note that $L_1$ and $L_2$ are cyclic difference sets 
 in $\F_{p^f}^\ast/C_0^{(k,p^f)}$ since $\chi_k(L_i)\overline{\chi_k(L_i)}=G_f(\chi_k)
\overline{G_f(\chi_k)}/p^{2s\theta}=p^{s(f-2\theta)}$. As determined in \cite{EHKX,SW02},  the corresponding cyclic $(v,k,\lambda)$ difference sets with 
$k\le (v-1)/2$ to cyclotomic 
strongly regular graphs of the Schmidt-White conjecture are as follows: 
\begin{enumerate}
\item[(1)] (subfield case) the Singer $(\frac{p^f-1}{p^d-1},\frac{p^{f-d}-1}{p^d-1},\frac{p^{f-2d}-1}{p^d-1})$ difference set;   
\item[(2)] (semi-primitive case) the trivial $(v,1,0)$ difference set;
\item[(3)] (exceptional case) see Table~\ref{T6}.
\begin{table}[h]
\caption{Cyclic difference sets corresponding to eleven sporadic examples}
\label{T6}
$$
\begin{array}{|c||c|c|c|c|}
\hline
\mbox{No.}&v&k&\lambda&\mbox{Name}\\
\hline
1&11&5&2&\mbox{Quadratic residue difference set~\cite[Theorem 1.12]{Beth}}\\
2&19&9&4&\mbox{Quadratic residue difference set}\\
3&35&17&8&\mbox{Twin-prime difference set~\cite[Theorem 8.2]{Beth}}\\
4&37&9&2&\mbox{Biquadratic residue difference set~\cite[Theorem 8.11]{Beth}}\\
5&43&21&10&\mbox{Hall's sextic difference set~\cite[Theorem 8.3]{Beth}}\\
6&67&33&16&\mbox{Quadratic residue difference set}\\
7&107&53&26&\mbox{Quadratic residue difference set}\\
8&133&33&8&\mbox{Quadratic residue difference set}\\
9&163&81&40&\mbox{Hall's sporadic difference set~\cite[Remarks 8.21(b)]{Beth}}\\
10&323&161&80&\mbox{Twin-prime difference set}\\
11&499&249&124&\mbox{Quadratic residue difference set}\\
\hline
\end{array}
$$
\end{table}
\end{enumerate}}
\end{remark}
The following is our main theorem of this paper. 
\begin{theorem} \label{thm:difmain1}
We assume that $ L_i\cap C_0^{(k,q)}=\emptyset$. 
Let $L'=\{y\,(\mod{k})\,|\,\gamma^{-y}\in L_i;0\le y\le q-2\}$ and 
let $I$ be the $|L_i|$-element set of odd integers modulo $2k$ such that  
$I\,(\mod{k})=L'$.  
Set 
\[
J=\{0\}\cup I \cup 
 2\left((\Z/k\Z)\setminus 2^{-1}\cdot (L'\cup \{0\})\right) \, (\mod{2k}). 
\]
Then, 
$D=\bigcup_{j\in J}C_j^{(2k,q)}$ in $\F_q$
 is a skew Hadamard  difference set 
or a Paley type partial difference set  according to $q\equiv 3\,(\mod{4})$ or 
$\equiv 1\,(\mod{4})$, i.e., it holds that 
\[
\psi(D)=\frac{-1\pm \sqrt{\pm q}}{2}.
\]  
\end{theorem}
\proof
First of all, we observe the following facts:  
\begin{description}
\item[(1)] It is clear that $J\,(\mod{k})=\{0,1,\ldots,k-1\}$. 
In particular,  if $q\equiv 3\,(\mod{4})$, i.e.,  $-1\in C_{k}^{(2k,q)}$, it follows that  $\F_q=\{0\}\cup D\cup -D$. 
\item[(2)] By the Davenport-Hasse product formula, it holds that 
\[
G_f(\chi_{2k})=\frac{G_f(\chi_{k})G_f(\chi_2)}{\chi_{k}(2)G_f(\chi_{k}^{2^{-1}})}. 
\]
Then, by noting that 
$G_f(\chi_{k}^{2^{-1}})G_f(\chi_{k}^{-2^{-1}})=\chi_{k}^{2^{-1}}(-1)q$ and 
the restriction of $\chi_{k}$ to $\F_p$ is trivial, it follows that
\begin{equation}\label{eq:Gauss2p1}
G_f(\chi_{2k})=\frac{1}{q}G_f(\chi_2)G_f(\chi_{k})G_f(\chi_{k}^{-2^{-1}}). 
\end{equation}
\item[(3)] The sum $\sum_{y\in J}\chi_{2k}^{x}(\gamma^y)$ for 
any $x$ such that $2,k\not |x$ is computable by 
using  Proposition~\ref{prop} as follows: 
\begin{eqnarray}
\sum_{y\in J}\chi_{2k}^{x}(\gamma^y)&=&
\sum_{y\in J}(-1)^y\chi_{k}^{x 2^{-1}}(\gamma^y)\nonumber\\ 
&=&1-\sum_{y\in L'}\chi_{k}^{x 2^{-1}}(\gamma^y)
+
\sum_{y\in (\Z/k\Z)\setminus (L'\cup \{0\})}\chi_{k}^{x 2^{-1}}(\gamma^y)\nonumber \\
&=&-2\sum_{y\in L'}\chi_{k}^{x2^{-1}}(\gamma^y)\nonumber\\
&=&-2\sum_{\omega\in L_i}\chi_{k}^{-x2^{-1}}(\omega)\nonumber\\
&=&(-1)^{i}2\epsilon G(\chi_k^{-x2^{-1}})/p^{s\theta}\label{eq:partialsum2}. 
\end{eqnarray}
\end{description}
Now, we compute the sum 
\[
T_a=\sum_{2,k\not |x}G_f(\chi_{2k}^{-x})
\sum_{y\in J}\chi_{2k}^{x}(\gamma^{a+y}). 
\]
By (\ref{eq:Gauss2p1}) and (\ref{eq:partialsum2}), we have 
\begin{eqnarray*}
T_a
&=&(-1)^{a+i}\epsilon\frac{2}{p^{s\theta}}\sum_{2,k\not |x}G_f(\chi_{2k}^{-x})G(\chi_k^{-x2^{-1}})\chi_{k}^{x2^{-1}}(\gamma^{a})\\
&=&(-1)^{a+i}\epsilon 2\frac{G_f(\chi_2)}{qp^{s\theta}}\sum_{x=1}^{k-1}
G_f(\chi_{k}^{-x})G_f(\chi_{k}^{x2^{-1}})
G(\chi_k^{-x2^{-1}})
\chi_{k}^{x2^{-1}}(\gamma^{a})\nonumber\\
&=&(-1)^{a+i} \epsilon 2\frac{G_f(\chi_2)}{p^{s\theta}}\sum_{x=1}^{k-1}
G_f(\chi_{k}^{-x})
\chi_{k}^{x}(\gamma^{2^{-1}a})\\
&=&(-1)^{a+i} \epsilon 2\frac{G_f(\chi_2)}{p^{s\theta}}
(k\cdot \psi(\gamma^{2^{-1}a} C_0^{(k,q)})+1),
\end{eqnarray*}
where we used  $G_f(\chi_{k}^{x2^{-1}})G_f(\chi_{k}^{-x2^{-1}})
=\chi_{k}^{x2^{-1}}(-1)q$.  
Then, by (\ref{eq:Wu1}), we obtain  
\begin{eqnarray*}
k(2\cdot \psi(aD)+1)&=&
\sum_{\ell=1}^{2k-1}G_f(\chi_{2k}^{-\ell})\sum_{y\in J}
\chi_{2k}^\ell(\gamma^{a+y})\\
&=&G_f(\chi_2)\sum_{y\in J}
\chi_{2}(\gamma^{a+y})+T_a\\
&=&(-1)^aG_f(\chi_2)\left(k-2|L_i|\right)\\
& &\hspace{1cm}+\epsilon 2(-1)^{a+i}\frac{G_f(\chi_2)}{p^{s\theta}}(k\cdot \psi(\gamma^{2^{-1}a} C_0^{(k,q)})+1)\\
&=&\pm (-1)^a k G_f(\chi_2). 
\end{eqnarray*}
By (\ref{eq:quad}), we obtain 
\[
\psi(aD)=\frac{-1\pm \sqrt{\delta q}}{2}, 
\]
where $\delta=1$ or $-1$ according to $q\equiv 3\,(\mod{4})$ or 
$\equiv 1\,(\mod{4})$. 
This completes the proof. \qed

Applying our theorem to subfield examples of 
cyclotomic strongly regular graphs, we obtain skew Hadamard 
difference sets in $\F_q$ for any $q=p^{st}$ with $st\ge 3$ by a nontrivial and different cyclotomic construction 
from that of the Paley difference sets although we do not know the 
constructed difference sets are inequivalent or not.   

Also, we may obtain an infinite family of skew Hadamard difference sets starting from each skew Hadamard difference  set of  Theorem~\ref{thm:difmain1} by applying the following theorem.  
\begin{theorem}(\cite{M12})\label{cor:sHd1}
Let $h=2p_1$ with an odd prime $p_1$ and let $p$ be a prime such that $\langle p\rangle$ is of index $e$ modulo $h$. Furthermore, 
let $k=2p_1^{m}$ and assume that $\langle p\rangle$ is again of index $e$ modulo $k$. 
Put $q=p^{(p_1-1)/e}$ and $q'=p^{p_1^{m-1}(p_1-1)/e}$. 
Define $J$ as any subset of $\{0,1,\ldots,h-1\}$ such that 
$J\,(\mod{p_1})=\{0,1,\ldots,p_1-1\}$. 
Let \[
D=\bigcup_{i\in J}C_{i}^{(h,q)} \mbox{\, and \, } 
D'=\bigcup_{i_1=0}^{p_1^{m-1}}\bigcup_{i\in J}C_{2i_1+ik/h}^{(k,q')}. 
\] 
If $D$ is a skew Hadamard difference set or a Paley type partial difference set  in $\F_{q}$, then  so does 
$D'$  in $\F_{q'}$.
\end{theorem}
By combining Theorems~\ref{thm:difmain1} and \ref{cor:sHd1}, we immediately have the following corollary, which yields an infinite family of skew Hadamard difference 
sets from  a cyclotomic strongly regular graph. 
\begin{corollary}\label{maincor}
Let $k=p_1^m$ and let $p$ be of index $e$ both of modulo $p_1$ and $k$. 
Put $q=p^{(p_1-1)/e}$, $q'=p^{p_1^{m-1}(p_1-1)/e}$, and 
\[
D=\bigcup_{i=0}^{p_1^{m-1}-1}\bigcup_{j\in J}C_{2i+p_1^{m-1}j}^{(2k,q)},
\]
where $J$ is defined as in Theorem~\ref{thm:difmain1}. 
If $\Cay(\F_q,C_0^{(p_1,q)})$  is strongly regular, then 
$D$ in $\F_{q'}$ is a skew Hadamard difference set 
or a Paley type partial difference set  according to $q\equiv 3\,(\mod{4})$ or 
$\equiv 1\,(\mod{4})$. 
\end{corollary}
\begin{example}\label{exam}{\em 
By Corollary~\ref{maincor}, we obtain new constructions of infinite families 
of skew Hadamard difference sets and Paley type partial difference sets 
for the quadruples $(p_1,p,f,e)$ of No. 2, 4, 5, 6, 7, 9, and 11 in Table~\ref{Tab1}. 
Note that  we can not obtain an infinite family of 
skew Hadamard difference sets from the  cyclotomic srg  of No. 1 because 
$p$ is not of index $2$ in $\Z/2p_1^m\Z$ for $m\ge 2$ while 
$\bigcup_{j\in \{0\}\cup \langle p\rangle \cup 2\langle p\rangle}C_j^{(2p_1,p^f)}$ forms a skew Hadamard difference set. 

Also, there are a lot of subfield examples satisfying 
$[(\Z/p_1\Z)^\ast:\langle p\rangle]=e$ and 
$p_1=\frac{p^{(p_1-1)/e}-1}{p^t-1}$ for some $t\,|\,(p_1-1)/e$. We list ten examples 
satisfying these conditions in Table~\ref{Tab4}. \begin{table}[h]
\caption{Subfield examples}
\label{Tab4}
$$
\begin{array}{|c|c|c|c||c|c|c|c|}
\hline
p_1&p&f&e&p_1&p&f&e\\
\hline
13&3&3&4&1723&41&3&574\\
31&5&3&10&2801&7&5&560\\
307&17&3&102&3541&59&3&1180\\
757&3&9&84&5113&71&3&1704\\
1093&3&7&156&8011&89&3&2670\\
\hline
\end{array}
$$
\end{table}
From these examples, we obtain  infinite families 
of skew Hadamard difference sets and Paley type partial difference sets by Corollary~\ref{maincor}.  
}
\end{example}
\section{Concluding remarks and open problems}
In this section, we give important remarks and open problems related to our results. 
\begin{remark}{\em 
In \cite{M12}, the author found two examples of skew Hadamard difference sets of index $4$, those are, 
$\bigcup_{j\in \{p_1\}\cup Q \cup 2Q}C_{j}^{(2p_1,p^f)}$  for 
$(p_1,p,f)=(13,3,3)$ and  $\bigcup_{j\in \{0\}\cup Q \cup 2Q}C_{j}^{(2p_1,p^f)}$   for 
$(p_1,p,f)=(29,7,7)$, 
where $Q$ is the subgroup of index $2$ of $(\Z/2p_1\Z)^\ast$. 
These two examples are not covered by Theorem~\ref{thm:difmain1}, i.e., 
there do not exist corresponding cyclotomic strongly regular graphs and cyclic difference sets. 
More generally, via a computation similar to \cite{GXY11}  involving known evaluations of index $4$ Gauss sums, one can prove that either of $\bigcup_{j\in \{0\}\cup Q \cup 2Q}C_{j}^{(2p_1,p^{(p_1-1)/4})}$ or 
 $\bigcup_{j\in \{p_1\}\cup Q \cup 2Q}C_{j}^{(2p_1,p^{(p_1-1)/4})}$  
is a skew Hadamard difference set or 
a Paley type  partial difference set 
in $\F_{p^{(p_1-1)/4}}$ if the following conditions are fulfilled: 
\begin{itemize}
\item[(i)] $p$ is of index $4$ modulo $p_1$, 
\item[(ii)] $p_1=4p^{(p_1-1)/4-2b}+1$, where $b$ is defined as 
\[
b=\min\left\{\frac{1}{p_1}\sum_{x\in S}x\,|\,S\in (\Z/p_1\Z)^\ast/\langle p\rangle \right\}, 
\] 
\item[(iii)] $p_1=A^2+4$ for some integer  
$A\equiv 3\,(\mod{4})$. 
\end{itemize}
The author found only three examples satisfying these conditions, which are
\[
(p_1,p,f)=(13,3,3),(29,7,7),(53,13,13). 
\] 
For each of these three examples, we obtain an infinite family 
of skew Hadamard difference sets or Paley type partial difference sets by applying 
Theorem~\ref{cor:sHd1}. Here, we have the following natural question. 
\begin{problem}
Determine for which $(p,p_1,e)$ either 
$\bigcup_{j\in \{0\}\cup Q \cup 2Q}C_{j}^{(2p_1,p^{(p_1-1)/e})}$ or 
 $\bigcup_{j\in \{p_1\}\cup Q \cup 2Q}C_{j}^{(2p_1,p^{(p_1-1)/e})}$  
forms a skew Hadamard difference set or a Paley type partial difference set. 
\end{problem}
Also, 
by computer, the author found an interesting example of a skew Hadamard 
difference set in  the case where $(p,f,p_1)=(7,3,19)$ and $e=6$:  
\[
D=\bigcup_{x\in I}C_{i}^{(2p_1,p^f)}, 
\]
where 
\[
I=\{p_1\}\cup \langle p\rangle\cup  3\langle p\rangle\cup 
3^3\langle p\rangle\cup 2\cdot 3\langle p\rangle\cup 2\cdot 3^3\langle p\rangle\cup 2\cdot 3^4\langle p\rangle\, (\mod{2p_1}).
\]
One can use a
computer to find that the automorphism group of the symmetric design $Dev(D)$ derived from $D$ has size $3^4\cdot 7^3$. (We will write the size as $\#Aut(Dev(D))$.) On the other hand, $\#Aut(Dev(P))=3^3\cdot 7^3\cdot 19$ for the Paley difference set $P$ with the same parameter. Thus, the skew Hadamard difference set $D$ is inequivalent to the Paley difference set.  Furthermore, 
since the size of the Sylow
$p$-subgroup of the automorphism group of the design derived from a
difference set constructed by Muzychuk~\cite{M} is strictly greater than $q$, we
conclude that $D$ is also inequivalent to the corresponding skew
Hadamard difference sets of \cite{M}.  Also, 
since 
the set $I$ satisfies $I\,(\mod{p_1})=\{0,1,\ldots,p_1-1\}$, we obtain an 
infinite family of skew Hadamard difference sets including this example by 
Theorem~\ref{cor:sHd1}. }
\end{remark}
\begin{remark}{\em 
As described in Introduction, to check  whether  obtained skew Hadamard difference sets and Paley type partial difference sets are equivalent or not to the classical Paley (partial) difference sets is very important. Although the problem is in general  difficult and the author could not prove that our construction always yields inequivalent skew Hadamard difference sets and Paley type partial difference sets to the Paley (partial) difference sets,  the author still believes that our infinite families include inequivalent ones abundantly. As an evidence for my believe, 
we can see by computer that  
the skew Hadamard 
difference set $D=\bigcup_{x\in J}C_{i}^{(2p_1,p^f)}$ 
with 
\[
J=\{0\}\cup \left(\bigcup_{i\in I}g^{i}\langle p\rangle\right)\cup 
 \left(2\bigcup_{i\in (\Z/e\Z)\setminus I}g^{i-s}\langle p\rangle\right)\hspace{0.5cm}(\mod{2p_1}), 
\]
where assume that $(\Z/k\Z)^\ast/\langle p\rangle $ is a cyclic group of order $e$ and let 
$g$ be a representative of a generator of $(\Z/k\Z)^\ast/\langle p\rangle $, 
is inequivalent to the Paley difference set in the following cases: 
\begin{itemize}
\item $(p,f,p_1)=(3,5,11)$, $(g,s)=(-1,1)$ and $I=\{0\}$: 
 In this case,  $\# Aut(Dev(D))=3^5\cdot 5\cdot 11$ and $\# Aut(Dev(P))=3^5\cdot 5\cdot 11^2$ for the corresponding Paley difference set $P$. 
\item $(p,f,p_1)=(3,7,1093)$, $(g,s)=(5,63)$ and take $I$ as $\bigcup_{i\in I}g^i \langle p\rangle =5\cdot (S+948))$ for the Singer difference set $S$ of PG$(6,3)$: In this case,  $\# Aut(Dev(D))=3^7\cdot 7$ and $\# Aut(Dev(P))=3^7\cdot 7\cdot 1093$ for the corresponding Paley difference set $P$. 
\item $(p,f,p_1)=(7,5,2801)$, $(g,s)=(3,58)$ and take $I$ as $\bigcup_{i\in I}g^i \langle p\rangle =3^{58}\cdot (S+292))$ for the Singer difference set $S$ of PG$(4,7)$: In this case,  $\# Aut(Dev(D))=3\cdot 5\cdot 7^5$ and $\# Aut(Dev(P))=3\cdot 5\cdot 7^5\cdot 2801$ for the corresponding Paley difference set $P$. 
\end{itemize} 

Furthermore, the author checked by computer that the Paley type srgs  with parameters $(p_1,p,f)=(31,5,3)$ and 
$(307,17,3)$ of  Example~\ref{exam}  are not isomorphic (as graph isomorphism) to the classical Paley graphs. 
(Note that   
in these cases  there is no factor $m>2$ of $p^f-1$ such that $p$ is semi-primitive modulo $m$.)  
\begin{problem}
Determine whether or not skew Hadamard difference sets and Paley type partial difference sets obtained in this paper  
are equivalent to the classical Paley (partial) difference sets. 
\end{problem}  }
\end{remark}


\section*{Acknowledgments}
The author would like to thank Tao Feng, Zhejian University, for his helpful 
comments on computations of the automorphism groups of 
symmetric designs derived from our skew Hadamard difference sets.   


\begin{thebibliography}{100}

\bibitem{BM73}
L. D. Baumert, J.~Mykkeltveit, Weight distributions of some irreducible 
cyclic codes, {\it DSN Progr. Rep.,} {\bf 16} (1973), 128--131.
 
\bibitem{BMW82}
L. D. Baumert, W. H. Mills, R. L. Ward, Uniform cyclotomy, {\it J. Number Theory}, {\bf 14} (1982), 67--82.

\bibitem{BEW97}
B. Berndt, R. Evans, K. S. Williams, {\it Gauss and Jacobi Sums}, Wiley, 1997. 
\bibitem{bjl} T. Beth, D. Jungnickel, H. Lenz, {\it Design Theory, vol. I, $2$nd  ed.,} Encyclopedia of Mathematics and its Applications 78, Cambridge University Press, Cambridge, 1999.
\bibitem{Beth}
T. Beth, D. Jungnickel, H. Lenz, {\it Design Theory}, Chapter VI, Cambridge University
Press, 1999.

\bibitem{bh} A. E. Brouwer,  W. H. Haemers, {\it Spectra of Graphs}, course notes, available at 
{\tt http://homepages.cwi.nl/{\~{ }}aeb/math/ipm.pdf }

\bibitem{CP}
Y. Q. Chen, J. Polhill, Paley type group schemes and planar Dembowski-Ostrom
polynomials, {\it Discr. Math.,} {\bf 311} (2011), 1349--1364.

\bibitem{DY06}
C. Ding, J. Yuan, A family of skew Hadamard difference sets, {\it J. Combin. Theory, Ser. A}, {\bf 113} (2006), 1526--1535.
\bibitem{DWX07} 
C. Ding, Z. Wang, Q. Xiang, Skew Hadamard difference sets from the Ree-Tits slice symplectic spreads in PG$(3, 3^{2h+1})$, {\it J. Combin. Theory, Ser. A}, {\bf 114} (2007), 867--887. 
\bibitem{EHKX}
R. Evans, H. Hollmann, C. Krattenthaler, Q. Xiang, Gauss sums, Jacobi sums, 
and $p$-ranks of cyclic difference sets,  {\it J. Combin. Theory, Ser. A}, {\bf 87} (1999), 74--119.
\bibitem{F11}
T. Feng, Non-abelian skew Hadamard difference sets fixed by a prescribed automorphism, {\it J. Combin. Theory, Ser. A}, {\bf 118} (2011), 27--36.
\bibitem{FX111}
T. Feng, Q. Xiang, Strongly regular graphs from union of cyclotomic classes,  to appear in {\it J. Combin. Theory, Ser. B}. 
\bibitem{FX113}
T. Feng, Q. Xiang, Cyclotomic constructions of skew Hadamard difference sets, {\it J. Combin. Theory, Ser. A}, {\bf 119} (2012), 245--256. 
\bibitem{FMX11} 
T. Feng, K. Momihara, Q. Xiang, Constructions of strongly regular Cayley graphs and skew
Hadamard difference sets from cyclotomic classes, 
ArXiv: 1201.0701. 
\bibitem{GXY11}
G. Ge, Q. Xiang, T. Yuan, Construction of strongly regular Cayley graphs 
using index four Gauss sums, ArXiv: 1201.0702. 


\bibitem{M94}
S. L. Ma, A survey of partial difference sets, {\it Des. Codes Cryptogr.,} {\bf 4} (1994), 221--261. 
\bibitem{M75}
R. J. McEliece,  Irreducible cyclic codes and Gauss sums,  in {\it Combinatorics}, pp. 183--200 ({\it Proc. NATO
Advanced Study Inst., Breukelen, 1974; M. Hall, Jr. and J. H. van Lint (Eds.)}), Part 1, Math. Centre Tracts, Vol.
55, Math. Centrum, Amsterdam, 1974. Republished by Reidel, Dordrecht, 1975 (pp. 185--202).
\bibitem{M}
M. E. Muzychuk, On skew Hadamard difference sets, ArXiv:1012.2089. 
\bibitem{M12} 
K. Momihara, 
Strongly regular Cayley graphs, skew Hadamard difference sets, 
and rationality of relative Gauss sums, ArXiv:1202.6414.
\bibitem{SW02}
B. Schmidt, C. White, All two-weight irreducible cyclic codes?, 
{\it Finite Fields Appl.,} {\bf 8} (2002), 321--367. 
\bibitem{S67}
T. Storer, {\it Cyclotomy and Difference Sets}, Lectures in Advanced Mathematics,  Markham Publishing Company, 1967.
\bibitem{VLSch} J. H. van Lint, A. Schrijver, Construction of strongly regular graphs, two-weight codes and partial geometries by finite fields, {\it Combinatorica,} {\bf 1} (1981), 63--73.
\bibitem{WH09}
G. B. Weng, L. Hu, Some results on skew Hadamard difference sets, {\it Des. Codes Cryptogr.,} {\bf 50} (2009), 93--105.
\bibitem{WQWX}
G. B. Weng, W. S. Qiu, Z. Wang, Q. Xiang, Pseudo-Paley graphs and skew
Hadamard difference sets from presemifields, {\it Des. Codes Cryptogr.,} 
{\bf 44} (2007),
49--62.
\bibitem{Wu12}
F. Wu, 
Constructions of strongly regular graphs using even
index Gauss sums, preprint.
\bibitem{YX10}
J. Yang, L. Xia, Complete solving of explicit evaluation of Gauss sums in the 
index $2$ case, {\it Sci. China Ser. A}, {\bf 53}  (2010), 2525--2542. 
\end{thebibliography}
\end{document}